%% file: ns.tex
\documentclass{article}

\newtheorem{theorem}{Theorem}
\newtheorem{definition}{Definition}
\newtheorem{lemma}[theorem]{Lemma}
\newtheorem{rem}{Remark}

\newfont{\bbc}{msbm10}
\def\Bbb#1{\hbox{{\bbc #1}}}
\def\qed{{\hfill{\vrule height5pt width3pt depth0pt}\medskip}}

\begin{document}
\begin{center}
\LARGE{\bf  Trapping regions  and an ODE-type proof of the
existence and uniqueness theorem for Navier-Stokes equations with
periodic boundary conditions on the plane }

\vskip 0.5cm

 Piotr Zgliczynski\footnote{Research supported in part by Polish KBN
grant 2P03A 011 18. %grant z Mrozkiem
} \\

\vskip 0.5cm

\today

\end{center}
\begin{center}
Jagiellonian University, Institute of Mathematics, \\ Reymonta 4,
30-059 Krak\'ow, Poland \\ e-mail: zgliczyn@im.uj.edu.pl \\ AND \\
Indiana University, Department of Mathematics, \\ 831 E Third
Street, Bloomington, IN 47405, USA \\ e-mail: pzgliczy@indiana.edu

\end{center}

\begin{abstract}
Using ODE-methods and trapping regions derived by Mattingly and Sinai
we give a new proof of  the existence and uniqueness of solutions
to Navier-Stokes equations with periodic boundary conditions on
the plane.
\end{abstract}



\section{Introduction.}

The goal of this paper to present self-contained account of the
ODE-type proofs from \cite{ES,MS,S} of the existence and
uniqueness of the Navier-Stokes systems with periodic boundary
conditions on the plane. Mattingly and Sinai  called their proof
elementary (see title of \cite{MS}), but  their proof was ODE-type
(elementary in their sense) only up to the moment of getting the
trapping regions for all Galerkin projections, but to pass to the
limit with the dimensions of Galerkin projections they invoked the
now standard results from \cite{CF,DG,T}  (which are not
elementary in any sense), which are usually not mastered by the
researchers working in dynamics of ODE's, to which this note is
mainly addressed. Here we fill-in this gap by giving ODE-type
arguments, which allow to pass to the limit. Using ODE-type
estimates based on the logarithmic norms we also obtained
uniqueness and an estimate for the Lipschitz constant of evolution
induced by the Navier-Stokes equations . In fact we have proved
that on the trapping region we have semidynamical system. The
results we prove here are well known for Navier-Stokes system in
2D (see for example \cite{FT,ES,K,DT}), but the method of getting
estimates on Galerkin projections presented in section
\ref{sec:limitlogn} appears to be new.

Another goal of this paper is to prepare the ground for the
rigorous study of the dynamics of the Navier-Stokes equations with
periodic boundary conditions. The trapping regions described here
are particular examples of the self-consistent apriori bounds
introduced in \cite{ZM} for the rigorous study of the dynamics of
the dissipative PDE's.

A few words about a general construction of the paper: In sections
\ref{sec:NSE} and \ref{sec:trap} we recall the results from
\cite{ES,MS,S} about the trapping regions. Sections
\ref{sec:limit} and \ref{sec:limitlogn} contain ODE-type proofs of
the convergence of the  Galerkin scheme on trapping regions. The
remaining sections contain the  existence results for the
Navier-Stokes equations on the plane and the Sannikov and Kaloshin
\cite{S} result in the dimension three.

\input feq.tex

\input ms.tex

\input limit.tex

\input limlogn.tex

\input exproof.tex

\input sankal.tex

\end{document}

%% file: feq.tex
\section{Navier-Stokes equations}
\label{sec:NSE}
 The general $d$-dimensional Navier-Stokes system
(NSS) is written for $d$ unknown functions
$u(t,x)=(u_1(t,x),\dots,u_d(t,x))$ of $d$ variables
$x=(x_1,\dots,x_d)$ and time $t$ and the pressure $p(t,x)$.
\begin{eqnarray}
  \frac{\partial u_i}{\partial t} + \sum_{k=1}^d u_k \frac{\partial u_i}{\partial
  x_k}= \nu \triangle u_i - \frac{\partial p}{\partial x_i} +
  f^{(i)} \label{eq:NS} \\
  \mbox{div}\ u= \sum_{i=1}^d \frac{\partial u_i}{\partial x_i}=0 \label{eq:div}
\end{eqnarray}
The functions $f^{(i)}$ are the components of the external
forcing, $\nu >0$ is the viscosity.

We consider (\ref{eq:NS}),(\ref{eq:div}) on the torus ${\Bbb T
}^d=\left({\Bbb R/2\pi}\right)^d$. This allows us to use Fourier
series. We write
\begin{equation}
  u(t,x)=\sum_{k \in {\Bbb Z}^d} u_k(t)e^{i(k,x)}, \qquad
  p(t,x)=\sum_{k \in {\Bbb Z}^d} p_k(t)e^{i(k,x)}
\end{equation}
Observe that $u_k(t) \in {\Bbb R}^d$, i.e. they are
$d$-dimensional vectors and $p_k(t) \in {\Bbb R}$. We will always
assume that $f_0=0$ and $u_0=0$.

Observe that (\ref{eq:div}) is reduced to the requirement $u_k
\bot k$. Namely
\begin{eqnarray*}
  \mbox{div}\ u= \sum_{k \in {\Bbb Z}^d} i (k,u_k(t))e^{i(k,x)} = 0 \\
    (k,u_k)=0 \quad k \in {\Bbb Z}^d
\end{eqnarray*}

To derive the evolution equation for $u_k(t)$ we will compute now
the nonlinear  term in (\ref{eq:NS}). We will use the following
notation
  $u_k=(u_{k,1},\dots,u_{k,d})$
\begin{eqnarray}
  \sum_l u_l \frac{\partial u}{\partial x_l}=\left(\sum_{k_1,l}
         u_{k_1,l}e^{i(k_1,x)}\right)\left( \sum_{k_2}i k_{2,l} u_{k_2}e^{i(k_2,x)}
         \right)= \\
=i \sum_{l,k_1,k_2}e^{i(k_1+k_2,x)} {k_2}_l \cdot u_{k_1,l} \cdot
u_{k_2} = i\sum_{k_1,k_2} e^{i(k_1+k_2,x)} (k_2|u_{k_1})u_{k_2} =
\\
i \sum_{k \in {\Bbb Z}^d}\left(\sum_{k_1}(u_{k_1}|k-k_1)u_{k-k_1}
\right)e^{i(k,x)}=i \sum_{k \in {\Bbb
Z}^d}\left(\sum_{k_1}(u_{k_1}|k)u_{k-k_1} \right)e^{i(k,x)}
\end{eqnarray}
We obtain the following infinite ladder of differential equations
for $u_k$
\begin{equation}
  \frac{d u_k}{d t}=-i \sum_{k_1}(u_{k_1}|k)u_{k-k_1} - \nu k^2u_k
  -i p_k k + f_k \label{eq:NSgal}
\end{equation}
Here $f_k$ are components of the external forcing. Let $\sqcap_k$
denote the operator of orthogonal projection to the
$(d-1)$-dimensional plane orthogonal to $k$. Observe that since
$(u_k,k)=0$ we have $\sqcap_k u_k=u_k$. We apply the projection
$\sqcap_k$ to (\ref{eq:NSgal}). The term $p_k k$ disappears and we
obtain
\begin{equation}
  \frac{d u_k}{d t}=-i \sum_{k_1}(u_{k_1}|k)\sqcap_k u_{k-k_1} - \nu k^2u_k + \sqcap_k f_k
    \label{eq:NSgal1}
\end{equation}
The pressure is given by the following formula
\begin{equation}
  -i \sum_{k_1}(u_{k_1}|k)(I - \sqcap_k)u_{k-k_1} -i p_k k + (I -
  \sqcap_k)f_k=0
\end{equation}

Observe that solutions of  (\ref{eq:NSgal1}) satisfy
incompressibility condition $(u_k,k)=0$. The the subspace of real
functions, which can be defined by $\overline{u_{-k}}=u_k$ for all
$k \in {\Bbb Z}^d$, where by $\overline{z}$ for $z \in {\Bbb C}$
we denote the conjugate of $z$, is invariant under
(\ref{eq:NSgal1}). In the sequel we will investigate the equation
(\ref{eq:NSgal1}) restricted to this subspace.

\begin{definition}
Enstrophy of $\{u_k, \quad k \in {\Bbb Z}^d\}$  is
\begin{displaymath}
  V(\{u_k, \quad k \in {\Bbb Z}^d\})= \sum_{k \in {\Bbb Z}^d}
  |k|^2|u_k|^2
\end{displaymath}
\end{definition}

%% file: ms.tex
\section{Construction of trapping regions from \cite{ES,MS}}
\label{sec:trap}

The idea in \cite{ES,MS} is to construct a trapping region for
each Galerkin projection and this trapping region give uniform
bounds allowing passing to the limit. The {\em trapping region }
for an ODE (here  Galerkin projection of Navier-Stokes equations)
is a set such that the vector field on its boundary is pointing
inside, hence no trajectory can leave it in forward time.
 In the sequel we consider
only the Galerkin projection onto the set of modes $O$, such that
if $ k \in O$ then $-k \in O$. We will call such projections {\em
symmetric}.

\begin{lemma}
$d=2$. For any solution of (\ref{eq:NSgal1}) (such that all
necessary Fourier  series converge) or the symmetric Galerkin
projection of (\ref{eq:NSgal1}) we have
\begin{equation}
  \frac{dV\{u_k(t)\}}{dt} \leq - 2\nu V(\{u_k(t) \}) + 2 V(F)\sqrt{ V(\{u_k(t) \})},
   \label{eq:enstrineq}
\end{equation}
where $V(F)=\sqrt{\sum |k|^2 f_k^2}$.
\end{lemma}

The proof can be found in many text-books, see also \cite{Si}.

The inequality (\ref{eq:enstrineq}) shows that
\begin{equation}
   \frac{dV\{u_k(t)\}}{dt} <0 , \qquad \mbox{when} \qquad
    V > V^*=\left(\frac{F}{ \nu} \right)^2 \label{eq:enstrineq2}
\end{equation}

\begin{lemma}
\label{lem:estmLin}
 Assume that $\{u_k, k \in{\Bbb Z}^d\}$ is such that for some $D <
 \infty$, $\gamma > 1 + \frac{d}{2}$
 \begin{equation}
   |u_k| \leq \frac{D}{k^\gamma}, \quad \mbox{and} \quad V(\{u_k\}) \leq
   V_0
 \end{equation}
 then for $d \geq 3$
 \begin{equation}
   | \sum_{k_1}(u_{k_1}|k)\sqcap_k u_{k-k_1}| \leq
   \frac{C \sqrt{V_0} D}{k^{\gamma - \frac{d}{2}}},
 \end{equation}
 where constant $C$ depends only on $\gamma$ and dimension $d$
 and for $d=2$ for any $\epsilon > 0$
 \begin{equation}
   | \sum_{k_1}(u_{k_1}|k)\sqcap_k u_{k-k_1}| \leq
   \frac{C(\epsilon,\gamma) \sqrt{V_0} D}{k^{\gamma - \frac{d}{2} - \epsilon}},
 \end{equation}
\end{lemma}
{\bf Proof:}

In order to estimate the sum $| \sum_{k_1}(u_{k_1}|k)\sqcap_k
u_{k-k_1}|$ we will use  the following inequality
\begin{equation}
    | (u_{k_1}|k)\sqcap_k u_{k-k_1}|=| (u_{k_1}|k - k_1)\sqcap_k u_{k-k_1}|
      \leq  |u_{k_1}| \ |k-k_1| \ |u_{k-k_1}| \\
\end{equation}

We consider three cases.

\noindent {\bf Case I.} $|k_1| \leq \frac{1}{2}|k|$.

Here $|k - k_1| \geq \frac{1}{2}|k|$ and therefore $|u_{k-k_1}|\
|k- k_1| \leq \frac{D}{|k-k_1|^{\gamma-1}} \leq \frac{2^{\gamma-1}
D}{|k|^{\gamma-1}}$. Now observe that
\begin{equation}
  \sum_{|k_1| \leq \frac{1}{2}|k|} |u_{k_1}| =
  \sum_{|k_1| \leq \frac{1}{2}|k|} |k_1| \ |u_{k_1}|
  \frac{1}{|k_1|} \leq \sqrt{\sum |k_1|^2 |u_{k_1}|^2} \cdot
  \sqrt{\sum_{|k_1| < \frac{1}{2}|k|} \frac{1}{|k_1|^2}}
\end{equation}

The sum $\sum_{|k_1| < \frac{1}{2}|k|} \frac{1}{|k_1|^2}$ can be
estimates from above by a constant times an integral  of
$\frac{1}{r^2}$ over the ball of radius $\frac{1}{2}|k|$ with the
ball around the origin removed. Hence for $d=2$ we have
\begin{equation}
  \sum_{|k_1| \leq \frac{1}{2}|k|} \frac{1}{|k_1|^2} \leq
  C\int_1^{|k|/2} \frac{r dr}{r^2} \leq C \ln|k|
\end{equation}

For $d \geq 3$ we have
\begin{equation}
  \sum_{|k_1| \leq \frac{1}{2}|k|} \frac{1}{|k_1|^2} \leq
  C\int_1^{|k|/2} \frac{r^{d-1} dr}{r^2} \leq C |k|^{d-2}
\end{equation}

From all the above computations it follows that for $d \geq 3$
holds
\begin{equation}
  | \sum_{|k_1| \leq \frac{|k|}{2}}(u_{k_1}|k)\sqcap_k u_{k-k_1}|
  \leq \frac{2^{\gamma -1} D}{|k|^{\gamma-1}} \sqrt{V_0} \sqrt{C}
  |k|^{\frac{d}{2} - 1} =
  \frac{2^{\gamma -1} D \sqrt{V_0} \sqrt{C} }{|k|^{\gamma - \frac{d}{2}}}
\end{equation}

For $d=2$ we have
\begin{equation}
  | \sum_{|k_1| \leq \frac{|k|}{2}}(u_{k_1}|k)\sqcap_k u_{k-k_1}|
  \leq \frac{2^{\gamma -1} D}{|k|^{\gamma-1}} \sqrt{V_0} \sqrt{C}
  \sqrt{ \ln |k|} < \frac{C \sqrt{V_0} D}{|k|^{\gamma -1 - \epsilon}}
\end{equation}

\noindent {\bf Case II. } $\frac{1}{2}|k| < |k_1| \leq 2 |k|$.
\begin{equation}
  |u_{k_1}| < \frac{D}{|k_1^\gamma|} < \frac{D}{\left(\frac{|k|}{2}\right)^\gamma}=
      \frac{2^\gamma D}{|k|^\gamma}
\end{equation}

Hence
\begin{equation}
  \sum_{\frac{1}{2}|k|  < |k_1| \leq 2 |k|}
     |u_{k_1}| \cdot |u_{k-k_1}|\cdot|k-k_1| \leq \frac{2^\gamma D}{|k|^\gamma}
  \sum_{\frac{1}{2}|k| < |k_1| \leq 2 |k|} |u_{k-k_1}|\cdot|k-k_1|
\end{equation}
We interpret $\sum_{\frac{1}{2}|k| < |k_1| \leq 2 |k|}
|u_{k-k_1}|\cdot|k-k_1|$ as a scalar product of
$|u_{k-k_1}|\cdot|k-k_1|$ and $1$, hence by the Schwarz inequality
\begin{equation}
\sum_{\frac{1}{2}|k| < |k_1| \leq 2 |k|} |u_{k-k_1}|\cdot|k-k_1|
\leq \sqrt{\sum_{|k_1| \leq 3|k|}|u_{k_1}|^2|k_1|^2} \cdot \sqrt{C
(3|k|)^d},
\end{equation}
where $C$ is such that $C(3|k|)^d$ is greater or equal than the
number of vectors in ${\Bbb Z}^d$, which are contained in the ball
of radius $3|k|$ around the origin.

Finally we obtain
\begin{equation}
 \sum_{\frac{1}{2}|k|  < |k_1| \leq 2 |k|}
     |u_{k_1}| \cdot |u_{k-k_1}|\cdot|k-k_1| \leq
    \frac{2^\gamma D {\tilde C} \sqrt{V_0}  }{|k|^{\gamma - \frac{d}{2}}}
\end{equation}

\noindent {\bf Case III.} $|k_1| > 2|k|$. Here we $|k - k_1| >
|k|$.
\begin{eqnarray*}
  \sum |u_{k_1}||k-k_1||u_{k-k_1}| \leq \frac{1}{|k|}\sum |u_{k_1}| |k_1| |k-k_1||u_{k-k_1}| \leq \\
   \frac{1}{|k|}\sqrt{\sum |u_{k_1}|^2|k_1|^2} \sqrt{\sum |u_{k-k_1}|^2 |k-k_1|^2} \leq
   \frac{\sqrt{V_0}}{|k|} \sqrt{\sum_{|k_1|>2|k|} \frac{D^2}{|k_1|^{2\gamma -2}} } = \\
   \frac{\sqrt{V_0} D}{|k|} \sqrt{\sum_{|k_1|>2|k|} \frac{1}{|k_1|^{2\gamma -2}} }
\end{eqnarray*}

To estimate $\sum_{|k_1|>2|k|} \frac{1}{|k_1|^{2\gamma -2}}$
observe that we have (we denote all constant factors depending on
$\gamma$ by $C$)
\begin{eqnarray*}
  \sum_{|k_1|>2|k|} \frac{1}{|k_1|^{2\gamma -2}} \leq
  C \int_{|k_1| > 2|k|} \frac{1}{|k_1|^{2\gamma-2}} d^d k_1 =
   C \int_{2|k|}^\infty \frac{1}{r^{2\gamma-2}}r^d dr = \\
   C \int_{2|k|}^\infty r^{-(2\gamma -2 -d +1)} = C |k|^{-2\gamma - 2 -d}
\end{eqnarray*}
Observe that we used here the assumption $\gamma > 1 +
\frac{d}{2}$, which guarantees that $2\gamma - 2-  d + 1 > 1 $ so
the integral converges.

Hence for the case III we obtain
\begin{equation}
  \sum ... \leq \frac{\sqrt{V_0} D C}{|k|^{\gamma - \frac{d}{2}}}
\end{equation}

Adding cases I,II,III we obtain for $d \geq 3$
\begin{equation}
 | \sum_{k_1}(u_{k_1}|k)\sqcap_k u_{k-k_1}| \leq
   \frac{C \sqrt{V_0}D}{ |k|^{ \gamma - \frac{d}{2} }  }
\end{equation}

For $d=2$ we obtain
\begin{equation}
 | \sum_{k_1}(u_{k_1}|k)\sqcap_k u_{k-k_1}| \leq
   \frac{C \sqrt{V_0}D}{ |k|^{ \gamma - \frac{d}{2} - \epsilon }  }
\end{equation}

\qed

\begin{lemma}
\label{lem:estmQ} Assume that $\gamma > d$, then
\begin{equation}
  \sum_{k_1 \in {\Bbb Z}^d \setminus \{0, k \}} \frac{1}{|k_1|^\gamma |k-k_1|^\gamma} \leq
    \frac{C_Q(d,\gamma)}{|k|^{\gamma}}
\end{equation}
\end{lemma}
{\bf Proof:}
 We consider three cases.

\noindent {\bf Case I.} $|k_1| < \frac{|k|}{2}$, hence $|k-k_1|
\geq \frac{|k|}{2}$.

We have
\begin{displaymath}
  \sum_{|k_1| < \frac{|k|}{2}} \leq \sum_{|k_1| < \frac{|k|}{2}}  \frac{1}{k_1^\gamma}
\frac{2^\gamma}{|k|^\gamma} <  \frac{2^\gamma}{|k|^\gamma} C
\int_1^{\infty}\frac{r^{d-1}}{r^\gamma} dr
\end{displaymath}
The improper integral $\int_1^{\infty} \frac{r^{d-1}}{r^\gamma}
dr$ converges, because $\gamma > d$.

Hence
\begin{displaymath}
 \sum_{|k_1| < \frac{|k|}{2}} < \frac{C_I(d,\gamma)}{|k|^\gamma}
\end{displaymath}

\noindent {\bf Case II.} $\frac{|k|}{2} < |k_1| \leq 2|k|$.
\begin{eqnarray*}
  \sum_{\frac{|k|}{2} < |k_1| \leq 2|k|} \leq \frac{2^\gamma}{|k|^\gamma}
  \sum_{\frac{|k|}{2} < |k_1| \leq 2|k|} \frac{1}{|k-k_1|^\gamma}<  \\
   \frac{2^\gamma}{|k|^\gamma}  \sum_{ |k_1| \leq 3|k|} \frac{1}{|k_1|^\gamma} <
   \frac{2^\gamma}{|k|^\gamma} C\int_1^{\infty}\frac{r^{d-1}}{r^\gamma} dr
\end{eqnarray*}

Hence
\begin{displaymath}
 \sum_{\frac{|k|}{2} < |k_1| \leq 2|k|} < \frac{C_{II}(d,\gamma)}{|k|^\gamma}
\end{displaymath}

\noindent {\bf Case III.} $2|k| < |k_1|$, hence $|k-k_1|>|k|$.
\begin{displaymath}
  \sum_{2|k| < |k_1|} < \frac{1}{|k|^\gamma} \sum \frac{1}{|k_1|^\gamma} <
  \frac{C_{III}(d,\gamma)}{|k|^\gamma}
\end{displaymath}

\qed

\subsection{The construction of the trapping region I.}
We take $V_0 > V^*$, $\gamma \geq 2.5$ and $K$ such that $f_k=0$
for $|k| > K$. We set
\begin{equation}
  N(V_0,K,\gamma,D)=\left\{  \{u_k\} \ | \ V(\{u_k\}) \leq V_0, \qquad
     |u_k| \leq \frac{D}{|k|^\gamma}, \quad |k| > K \right\}
\end{equation}

We prove that
\begin{theorem}
\label{thm:trap1} Let $d=2$ and $C(\gamma)$ be a constant from
lemma \ref{lem:estmLin}.
 If $K > \frac{C^2 V_0}{\nu^2}$ and $D
> \sqrt{V_0} K^{\gamma -1}$, then $N=N(V_0,K,\gamma,D)$ is a
trapping region for each Galerkin projection.
\end{theorem}
{\bf Proof:}  Observe that for $D \geq \sqrt{V_0} K^{\gamma -1}$
for all $\{u_k\} \in N$ holds
\begin{equation}
  |u_k| \leq \frac{D}{|k|^\gamma}. \label{eq:prdec}
\end{equation}
To prove this observe that (\ref{eq:prdec}) holds for $|k| > K$ by
the definition of $N$. For $|k| \leq K$ we proceed as follows:
since $V(\{u_k\}) \leq V_0$ then $|k|^2|u_k|^2 \leq V_0$. So we
have
\begin{equation}
  |u_k| \leq \frac{\sqrt{V_0}}{|k|} \leq \frac{D}{|k|^\gamma}, \quad |k| \leq K
\end{equation}
for $D$ such that $\sqrt{V_0}|k|^{\gamma-1} \leq D$ for all $|k|
\leq K$.

We will show now that on the boundary of $N$ (we are considering
the Galerkin projection) the vector field is pointing inside. For
points $V(\{u_k\})=V_0$ it follows from (\ref{eq:enstrineq2}). For
points such that $u_k=\frac{D}{|k|^\gamma}$ for some $|k| > K$ we
have from lemma \ref{lem:estmLin} (with $\epsilon=1/2$)
\begin{equation}
  \frac{d|u_k|}{dt} \leq
  \frac{C \sqrt{V_0} D}{|k|^{\gamma - \frac{3}{2}}} - \nu |k|^2 \frac{D}{|k|^{\gamma}} <0,
\end{equation}
which is satisfied when
\begin{equation}
   C \sqrt{V_0}  < \nu |k|^{1/2} . \label{eq:isobigK}
\end{equation}
Observe that (\ref{eq:isobigK}) holds for $|k| \geq K$ if $K >
\frac{C^2 V_0}{\nu^2}$.

\qed

\begin{rem}
Observe that in the proof it was of crucial importance that the
constant  $D$ entered linearly in the estimate in lemma
\ref{lem:estmLin} and due to this fact did not appear in
(\ref{eq:isobigK}). For example assume that the estimate of the
nonlinear part will be of the form $\frac{D^2 C }{|k|^{\gamma -
\frac{3}{2}}}$ then instead of (\ref{eq:isobigK}) we will have
\begin{displaymath}
   C D  < \nu |k|^{1/2}
\end{displaymath}
which will require that $K > \frac{C^2 D^2}{\nu^2}$, which might
be incompatible with $D > \sqrt{V_0}K^{\gamma-1}$.

This shows how important it was to use the enstrophy in these
estimates.
\end{rem}

\subsection{The construction of the trapping region II - exponential decay}
\begin{theorem}
\label{thm:trap2} Assume that $\gamma \geq 2.5$, $d=2$. Then the
set
\begin{equation}
  N_e=N(V_0,K,\gamma,D) \cap \left\{ \{ u_k \}\ | \ |u_k| \leq \frac{D_2}{|k|^\gamma}e^{-a|k|} \
  \mbox{for $|k| > K_e$} \right\},
\end{equation}
where $N(V_0,K,\gamma,D)$ is a trapping region from theorem
\ref{thm:trap1}, $D_2 > D$, $K_e
> \frac{C_Q(d,\gamma)D_2}{\nu}$ ($C_Q$ was obtained in lemma
\ref{lem:estmQ}) and $0 < a < \frac{1}{K_e}\ln \frac{D_2}{D}$ is a
trapping region for each Galerkin projection.
\end{theorem}
{\bf Proof:} The set $N_e$ constructed so that for all $|k|\leq
K_e$ the trapping (the vector field is pointing toward the
interior of $N_e$ on the boundary) is obtained from
$N(V_0,K,\gamma,D)$ and for $|k|
> K_e$ it results from the new exponential estimates.

Observe that $a$ is such that $\frac{D_2}{|k|^\gamma}e^{-a|k|} >
\frac{D}{|k|^\gamma}$ for all $|k| \leq K_e$. This solves the
trapping for $|k| \leq K_e$.

Hence to prove the trapping it is enough to consider the boundary
points such that $|u_k|=\frac{D_2}{|k|^\gamma}e^{-a|k|}$ for some
$k
> K_e$. For such a point and $|k|$ we have
\begin{eqnarray*}
  \frac{d |u_k|}{dt} \leq \left|\sum (u_{k_1}|k) \sqcap_k u_{k-k_1}\right| - \nu |k|^2 |u_k|
    \leq \\
    \sum |u_{k-1}||k||u_{|k-k_1|}| - \nu |k|^2 |u_k| \leq
    D_2^2 |k| \sum \frac{e^{-a|k_1|} e^{-a|k-k_1|}}{|k_1|^\gamma |k-k_1|^\gamma} -
     \nu |k|^2 |u_k|
\end{eqnarray*}
Observe that $e^{-a|k_1|} e^{-a|k-k_1|} \leq e^{-a|k|}$. From this
and lemma \ref{lem:estmQ} we obtain
\begin{displaymath}
  \frac{d |u_k|}{dt} < \frac{D_2^2 C_Q(\gamma,d)}{|k|^{\gamma-1}}e^{-a|k|} - \nu |k|^2 |u_k|
\end{displaymath}
Hence $ \frac{d |u_k|}{dt}  < 0$, when
\begin{displaymath}
|u_k|=\frac{D_2}{|k|^\gamma}e^{-a|k|} > \frac{C_Q D_2^2}{\nu
|k|^{\gamma+1}}e^{-a|k|}
\end{displaymath}

Which is equivalent to
\begin{displaymath}
 |k| > K_e=\frac{C_Q D_2}{\nu}.
\end{displaymath}

\qed

\subsection{Trapping region III - exponential decay in time }
\begin{theorem}
\label{thm:trap3}  Let $t_0 >0$. Assume that $\gamma \geq 2.5$,
$d=2$. Then the set
\begin{equation}
  N_e=N(V_0,K,\gamma,D) \cap \left\{ \{ u_k \}\ | \ |u_k| \leq \frac{D_3}{|k|^\gamma}e^{-a_3|k|t} \
  \mbox{for $|k| > K_e$} \right\},
\end{equation}
where $N(V_0,K,\gamma,D)$ is a trapping region from theorem
\ref{thm:trap1}, $D_3 > D$, $K_e> \frac{D_3 C_Q(d,\gamma)}{\nu} $
($C_Q$ was obtained in lemma \ref{lem:estmQ}) and $0 < a_3 <
\frac{1}{K_e t_0} \ln \frac{D_3}{D} $ is a trapping region for
each Galerkin projection for $0 \leq t \leq t_0$.
\end{theorem}
{\bf Proof:} The set $N_e$ constructed so that for all $|k|\leq
K_e$ the trapping is obtained from $N(V_0,K,\gamma,D)$ and for
$|k|
> K_e$ it results from the new exponential estimates.

To be sure that the boundary of $N_e$ for $|k| < K_e$ is obtained
from \newline $N(V_0,K,\gamma,D)$ we require that
\begin{equation}
  \frac{D}{|k|^\gamma} <  \frac{D_3}{|k|^\gamma}e^{-a_3|k|t},
   \quad \mbox{for $0 \leq t \leq t_0$ and $|k|  \leq K_e$}. \label{eq:tra31}
\end{equation}
Easy computations show that (\ref{eq:tra31}) holds iff $ a_3 <
\frac{1}{K_e t_0} \ln \frac{D_3}{D}$.

To have the trapping for $|k| > K_e$ we need to show that$\frac{d
|u_k|}{dt} < 0$ if $|u_k|=\frac{D_3}{|k|^\gamma}e^{-a_3 t}$, for
some $0 \leq t \leq t_0$ and $|k| > K_e)$.

\begin{eqnarray*}
\frac{d |u_k|}{dt} \leq \sum |u_{k_1}||k||u_{k-k_1}| - \nu |k|^2
|u_k| \leq \\ |k| D_3^2 \sum \frac{e^{-a_3|k_1|t}
e^{-a_3|k-k_1|t}}{|k_1|^\gamma |k-k_1|^\gamma} - \nu |k|^2 |u_k|
\leq \\ |k|e^{-a_3|k|t} D_3^2 \sum \frac{1}{|k_1|^\gamma
|k-k_1|^\gamma} -  \nu |k|^2 |u_k| \leq \\
\frac{e^{-a_3|k|t}
D_3^2C_Q(d,\gamma)}{|k|^{\gamma -1}} - \nu |k|^2 |u_k|
\end{eqnarray*}
Hence $\frac{d |u_k|}{dt} < 0$ if
\begin{equation}
  \frac{D_3^2 C_Q(d,\gamma)}{\nu |k|^{\gamma +1}}e^{-a_3|k|t} <
     |u_k|=\frac{D_3}{|k|^\gamma}e^{-a_3|k|t},
\end{equation}
which is equivalent to
\begin{equation}
  \frac{D_3 C_Q}{\nu} < |k|.
\end{equation}
Hence for $K_e \geq  \frac{D_3 C_Q}{\nu}$ we obtain the trapping.
 \qed

%% file: limit.tex
\section{Passing to the limit for Galerkin projections via Ascoli-Arzela lemma}
\label{sec:limit}

The goal of this section is a relatively simple argument for the
passing to the limit with Galerkin projections.

All what follows  was essentially proved in \cite{ZM}. We will
also use some conventions used there.

Let $H$ be a Hilbert space. Let $\phi_1,\phi_2,\dots$ be a
orthonormal basis in $H$.

Let $A_n:H \to H$ be denote the projection onto 1-dimensional
subspace $\langle \phi_n \rangle$, i.e $x = \sum A_n(x) \phi_n$
for all $x \in H$. By $V_n$ we will denote  the space spanned by
$\{\phi_1,\dots,\phi_n\}$.  Let $P_n$ denote the projection onto ,
$Q_n=I-P_n$.

\begin{definition}
Let  $W\subset H$ and $F:\mbox{dom}(H) \to H$. We say that $W$ and
$F$ satisfy conditions C1,C2,C3 if
\begin{description}
\item[C1] There exists $M \geq 0$, such that $P_n(W) \subset W$ for $k \geq M$
\item[C2] Let $\hat{u_k}=\max_{x \in W} |A_k x|$. Then, $\hat{u} = \sum \hat{u_k}\in H$.
  In particular, $|\hat{u}|<\infty$.
\item[C3] The function $x\mapsto F(x)$ is continuous on $W$. The sequence
   $f=\{f_k\}$, given by $f_k=\max_{x \in W} |A_k F(x)|$ is in $H$.
  In particular, $|f|<\infty$.
\end{description}
\end{definition}

Observe that  condition C2 implies that the set $W$ is compact.
Conditions C2 and C3 guarantee good behavior of $F$ with respect
to passage the limit. We have here continuous function on the
compact set, this is also perfect setting for study the dynamics
of $x'=F(x)$ (see \cite{ZM} for more details).

\begin{lemma}
\label{lem:ladder} Assume that  $W \subset H$ and $F$ satisfy
C1,C2,C3. Let $x:[0,T] \to W$ be such that for each $n$
\begin{equation}
   \frac{d A_n x}{dt}=A_n(F(x)).
\end{equation}
Then
\begin{equation}
  x'=F(x).
\end{equation}
\end{lemma}
{\bf Proof:} Let us set $x_k=A_k x$. Let us fix $\epsilon > 0$ and
$t \in [0,T]$. For any $n$ we have
\begin{eqnarray}
  \left|\frac{x(t+h) - x(t)}{h} - F(x)\right| \leq
    \left|\frac{P_n x(t+h) - P_n x(t)}{h} - P_n F(x)\right| + \\
      \left|\frac{1}{h} \sum_{k=n+1}^\infty (x_k(t+h) - x_k(t))\phi_k  \right| + |Q_n F(x)|
\end{eqnarray}
We will estimate the three terms on the right hand side
separately. From {\bf C3} it follows for a given $\epsilon
>0$ there exists $n_0$ such that $n>n_0$ implies
\[ \left| Q_n(F(a)) \right|  < \epsilon/3. \]
From now on fix $n>n_0$. Again {\bf C3} and the mean value theorem
implies
\begin{eqnarray*}
\left|\sum_{k=n+1}^{\infty} \frac{1}{h}(x_k(t+h)
-x_k(t))\phi_k\right| &=& \left|\sum_{k=n+1}^{\infty} \frac{d
x_k}{d t}(t+\theta_k h) \phi_k\right| \\ &\leq&
   \left| \sum_{k=n+1}^\infty f_k \phi_k \right|   < \epsilon/3 .
\end{eqnarray*}
Finally, for $h$ sufficiently small,
\[
\left|\frac{1}{h} (P_nx(t+h) -P_nx(t))\phi_k - P_n F(x) \right|
  < \epsilon/3
\]
and hence the desired limit is obtained.

 \qed

\begin{lemma}
\label{lem:limit} Assume that $W \subset H$ and the function $F$
satisfy C1,C2,C3. Let $x_0 \in W$. Assume that for each $n$ a
function $x_n:[0,T] \to P_n(W)$ is a solution of the problem
(Galerkin projection of $x'=F(x)$)
\begin{equation}
   x_n'=P_n(F(x)), \quad x_n(0)=P_n(x_0).
\end{equation}
Assume also that $x_n$ converges uniformly to $x^*:[0,T] \to W$.

Then $x^*$ solves the following initial value problem
\begin{equation}
  x'=F(x), \quad x(0)=x_0
\end{equation}
\end{lemma}
{\bf Proof:} We show first  that for all $n$ and $t \in [0,T]$
holds
\begin{equation}
  P_nx^*(t)=P_n x_0 + \int_{0}^{t} P_nF(x^*(s)) ds. \label{eq:int}
\end{equation}

Let us fix $n$. Observe that for each $m \geq n$ the following
equation holds
\begin{equation}
  P_n x_m(t)=P_nx_0 + \int_{0}^t P_n F(x_m(s)) ds \label{eq:intGals}
\end{equation}
Since the series $x_m$ converges uniformly to $x^*$, then also
$P_n x_m$ converges uniformly to $P_n x^*$. Observe that also the
functions $P_nF(x_m)$ converge uniformly to $P_n F(x^*)$ as the
composition of the uniformly continuous function $P_n F$ (because
F is a continuous function on the compact set $W$) with a
uniformly convergent sequence, hence also the integral in
(\ref{eq:intGals}) is converging (uniformly in $t \in [0,T]$) to $
\int_{0}^{t} P_nF(x^*(s))$. This proves (\ref{eq:int}).
Differentiation of (\ref{eq:int}) gives
\begin{equation}
  \frac{d P_n x^*}{dt}=P_nF(x^*).
\end{equation}
The assertion follows from lemma \ref{lem:ladder}  \qed

\begin{theorem}
\label{thm:limit} Assume that $W \subset H$ and the function $F$
satisfy C1,C2,C3. Let $x_0 \in W$. Assume that for each $n$ a
function $x_n:[0,T] \to P_n(W)$ is a solution of the problem
(Galerkin projection of $x'=F(x)$)
\begin{equation}
   x_n'=P_n(F(x)), \quad x_n(0)=P_n(x_0).
\end{equation}

Then there exists $x^*:[0,T] \to W$, such that $x^*$ solves the
following initial value problem
\begin{equation}
  x'=F(x), \quad x(0)=x_0 \label{eq:feq}
\end{equation}
\end{theorem}
{\bf Proof:} The idea goes as follows, we would like to pickup a
convergent subsequence from $\{x_n\}$ using Ascoli-Arzela
compactness theorem. Later we show that the limit function $x^*$
solves (\ref{eq:feq}).

Observe first that due to compactness of $W$ and since $x_n(t) \in
W$ for $t \in [0,T]$ the sequence $\{x_n\}$ is contained in a
compact set. Observe that the derivatives $x_n'(t)$ are  uniformly
bounded by $|F(W)|$, hence the sequence of functions $x_n$ is
equicontinuous. From Ascoli-Arzela theorem if follows that there
exists a subsequence converging uniformly to $x^*:[0,T]\to W$.
Without loss of generality we can assume that the whole sequence
$x_n$ is converging uniformly to $x^*$. It is obvious that
$x^*(0)=x_0$. The assertion of the theorem follows from lemma
\ref{lem:limit}. \qed

%% file: limlogn.tex
\section{Passing to the limit, an analytic argument}
\label{sec:limitlogn}

The goal of this section is to present another argument for the
limit of Galerkin projections. Compared to section \ref{sec:limit}
we assume more about the function $F$ and we add a new condition D
on the trapping regions, which are satisfied for the Navier-Stokes
system and the trapping regions constructed in section
\ref{sec:trap}. We obtain better results about convergence plus
uniqueness and Lipschitz constants for the induced flow.

 We will use here the notations introduced in
section \ref{sec:limit}. We investigate the Galerkin projections
of the following problem
\begin{equation}
  x'=F(x)=L(x) + N(x), \label{eq:feqLN}
\end{equation}
where $L$ is a linear operator and $N$ is a nonlinear part of F.
We assume that the basis $\phi_1, \phi_2, \dots$ of $H$ is build
from eigenvectors of $L$. We assume that the corresponding
eigenvalues $\lambda_k$ (i.e. $L\phi_k=\lambda_k \phi_k$) can
ordered so that
\begin{displaymath}
  \lambda_1 \geq \lambda_2 \geq \dots, \quad \mbox{and} \quad \lim_{k\to \infty} \lambda_k=-\infty.
\end{displaymath}
Hence we can have only a finite number of positive eigenvalues.

\subsection{Estimates based on logarithmic norms }

The goal of this section is to recall some results about one-sided
Lipschitz constants of the flows induced by ODE's. We will invoked
here results from \cite{HNW}.

\begin{definition} \cite[Def. I.10.4]{HNW}
Let $Q$ be a square matrix, then we call
\begin{displaymath}
  \mu(A) = \lim_{h > 0 ,h\to 0} \frac{\|I + hQ \| - 1}{h}
\end{displaymath}
the {\em logarithmic norm} of $Q$.
\end{definition}

\begin{theorem} \cite[Th. I.10.5]{HNW}
The logarithmic norm is obtained by the following formulas
\begin{itemize}
\item for Euclidean norm
  \begin{displaymath}
     \mu(Q)=\mbox{the largest eigenvalue of } \quad 1/2(Q + Q^T).
  \end{displaymath}
\item for max norm $\|x\|=\max_{k} |x_k|$
  \begin{displaymath}
     \mu(Q)=\max_k \left( q_{kk} + \sum_{i\neq k} |q_{ki}| \right)
  \end{displaymath}
\item for norm $\|x\|=\sum_k  |x_k|$
  \begin{displaymath}
     \mu(Q)=\max_i \left( q_{ii} + \sum_{k\neq i} |q_{ki}| \right)
  \end{displaymath}
\end{itemize}
\end{theorem}

Consider now the differential equation
\begin{equation}
  x'=f(x), \quad \mbox{$f \in C^1$}. \label{eq:odelogn}
\end{equation}
Let $\varphi(t,x_0)$ we denote the solution of equation
(\ref{eq:odelogn}) with the initial condition $x(0)=x_0$. By $\|x
\|$ we denote a fixed arbitrary norm in ${\Bbb R}^n$.

The following theorem was proved in \cite[Th. I.10.6]{HNW} (for
nonautonoumous ODE, here we restrict ourselves  to an autonomous
case only and we use a different notation)
\begin{theorem}
\label{th:lognorm} Let $y:[0,T] \to {\Bbb R}^n$  and
$\varphi(\cdot,x_0)$ is defined for $t \in [0,T]$. Suppose that we
have the following estimates
\begin{eqnarray*}
  \mu\left(\frac{\partial f}{\partial x}(\eta)\right) \leq l(t),\quad \mbox{ for $\eta \in [y(t),x(t)] $} \\
  \left\| \frac{dy}{dt}(t) - f(y(t)) \right\| \leq \delta(t), \quad \|y(0) - x_0 \| \leq \rho.
\end{eqnarray*}
Then for $0 \leq t \leq T$ we have
\begin{displaymath}
 \| \varphi(t,x_0) - y(t)  \| \leq e^{L(t)}\left( \rho + \int_0^te^{-L(s)}\delta(s)ds  \right),
\end{displaymath}
where $L(t)=\int_0^tl(s)ds$.
\end{theorem}

From the above theorem one easily derives the following
\begin{lemma}
\label{lem:estmLogN}
Let $y:[0,T] \to {\Bbb R}^n$  and
$\varphi(\cdot,x_0)$ is defined for $t \in [0,T]$. Suppose that
$Z$ is a convex set such that  we have the following estimates
\begin{eqnarray*}
  y([0,T]), \varphi([0,T], x_0) \in Z  \\
  \mu\left(\frac{\partial f}{\partial x}(\eta)\right) \leq l,\quad \mbox{ for $\eta \in Z$} \\
  \left\| \frac{dy}{dt}(t) - f(y(t)) \right\| \leq \delta, \quad \|y(0) - x_0 \| \leq \rho.
\end{eqnarray*}
Then for $0 \leq t \leq T$ we have
\begin{displaymath}
 \| \varphi(t,x_0) - y(t)  \| \leq e^{lt}\rho + \delta \frac{e^lt -1}{l},\quad \mbox{if $l \neq 0$}.
\end{displaymath}
For $l=0$ we have
\begin{displaymath}
 \| \varphi(t,x_0) - y(t)  \| \leq \rho + \delta t.
\end{displaymath}
\end{lemma}

\subsection{Application to Galerkin projections - uniqueness and another proof of convergence}
\begin{definition}
We say  that $W \subset H$  and $F=N+L$ satisfy condition {\em D}
if the following condition holds
\begin{description}
\item[D] there exists  $l \in {\Bbb R}$ such that for all $k=1,2,\dots$
\begin{equation}
  1/2\sum_{i=1}^\infty {\left|\frac{\partial N_k}{\partial x_i}\right|}(W) +
  1/2 \sum_{i=1}^\infty {\left|\frac{\partial N_i}{\partial x_k}\right|}(W)  + \lambda_k \leq l
\end{equation}
\end{description}
\end{definition}

The main idea behind the condition {\bf D} is the ensure that the
logarithmic norms for all Galerkin projections are uniformly
bounded.

\begin{theorem}
\label{thm:limitLN}
  Assume that $W \subset H$ and $F$ satisfy conditions {\em C1,C2,C3,D} and
  $W$ is convex. Assume that $P_n(W)$ is a trapping region for the $n$-dimensional
  Galerkin projection of (\ref{eq:feqLN}) for all $n>M_1$. Then
  \begin{description}
    \item[1.]{\bf Uniform convergence and existence} For a fixed $x_0 \in W$,
    let  $x_n:[0,\infty] \to P_n(W)$ be a solution of $x'=P_n(F(x))$, $x(0)=P_nx_0$.
    Then $x_n$ converges uniformly on compact  intervals to a function
     $x^*:[0,\infty]\to W$, which
    is a solution of (\ref{eq:feqLN}) and $x^*(0)=x_0$
    \item[2.] {\bf Uniqueness within $W$}. There exists only one solution of the initial
     value problem (\ref{eq:feqLN}), $x(0)=x_0$ for any $x_0 \in W$, such that $x(t) \in W$
     for $t>0$.
    \item[3.]  {\bf Lipschitz constant}. Let $x:[0,\infty] \to W$ and
      $y:[0,\infty] \to W$ be solutions of (\ref{eq:feqLN}), then
      \begin{displaymath}
         |y(t) - x(t) | \leq e^{lt}|x(0) - y(0)|
      \end{displaymath}
    \item[4.] {\bf Semidynamical system.} The map $\varphi:{\Bbb R}_+ \times W \to W$, where
       $\varphi(\cdot,x_0)$ is a unique solution of equation (\ref{eq:feqLN}), such that
       $\varphi(0,x_0)=x_0$ defines a semidynamical system on $W$, namely
        \begin{itemize}
          \item $\varphi$ is continuous
          \item $\varphi(0,x)=x$
          \item $\varphi(t,\varphi(s,x)) = \varphi(t+s,x)$
        \end{itemize}
  \end{description}
\end{theorem}
{\bf Proof:}
By $|x|_n$ we will denote $|P_n(x)|$, i.e. Euclidean
norm in ${\Bbb R}^n$.

Let
\begin{displaymath}
\delta_n= \max_{x \in W} |P_n(F(x)) - P_n(F(P_n x))|.
\end{displaymath}
 Obviously $\delta_n \to 0$ for $n \to \infty $, because $F$ is uniformly continuous
 on $W$ and $P_n W \subset W$, for $n \geq M$.

Let us consider the logarithmic norm of the vector field for the
$n$-dimensional Galerkin projection. We will estimate it using the
euclidean norm  on $P_nH={\Bbb R}^n$(which coincides with the norm
inherited from $H$). Since
\begin{equation}
  \left[\frac{\partial P_n (L+N)}{\partial( x_1\dots x_n)}\right]_{ij}=
  \frac{\partial N_i}{\partial x_j} + \delta_{ij} \lambda_j,
\end{equation}
we need to estimate the largest eigenvalue of the following
matrix, $Q_n(x)$ for $x \in P_n(W)$,
\begin{equation}
  Q_{n,ij}(x)=\frac{1}{2} \frac{\partial N_i}{\partial x_j}(x) +
    \frac{1}{2} \frac{\partial N_j}{\partial x_i}(x) +  \delta_{ij} \lambda_j,
   \quad \mbox{for $i,j=1,\dots,n$}
\end{equation}
where $\delta_{ij}$ is a Kronecker symbol, i.e. $\delta_{ij}=1$,
if $i=j$ and $\delta_{ij}=0$ otherwise.

To estimate the largest eigenvalue of $Q_n$ we will use  the
Gershgorin theorem (see \cite[Property 5.2]{QSS}), which states
that all  eigenvalues of  a square $n \times n$-matrix $A$,
$\sigma(A)$, satisfy
\begin{equation}
  \sigma(A) \subset \cup_{j=1}^n \{z \in {\Bbb C}: \quad
  |z-A_{jj}| <  \Sigma_{i=1,i\neq j} |A_{ij}| \}.
\end{equation}
From above equation and condition D it follows immediately that
eigenvalues of $Q_n$ are less then or equal to $l_n$, where
\begin{equation}
  l_n = \max_{k=1,\dots,n} \max_{x \in P_n W}  \sum_{i=1}^n \left(
   1/2 \left| \frac{\partial N_k}{\partial x_i}(x)  \right| +
   1/2 \left| \frac{\partial N_i}{\partial x_k}(x)  \right| \right)  + \lambda_k.
\end{equation}
From the assumption {\bf D}, it follows that we have uniform bound
on $l_n$, namely
\begin{equation}
  l_n \leq l, \quad \mbox{for all $n$}.
\end{equation}

Let us take $m \geq n$. Let $x_n:[0,T] \to P_nW$ and $x_m:[0,T]
\to P_m W$ be the solutions of $n$- and $m$-dimensional
projections of (\ref{eq:feqLN}). From lemma \ref{lem:estmLogN} it
follows immediately that (we treat here $P_n x_m$ as a perturbed
'solution' $y$)
\begin{equation}
  |x_n(t) - P_n(x_m(t))|_n \leq e^{lt}|x_n(0) - P_nx_m(0)| + \delta_n \frac{e^{lt} - 1}{l}
  \label{eq:estmdiffg}
\end{equation}

To prove uniform convergence of $\{x_n\}$ starting from the same
initial condition  observe that
\begin{eqnarray*}
 |x_n(t) - x_m(t)| \leq  |x_n(t) - P_n(x_m(t))|_n + |(I-P_n) x_m(t)| \leq \\
    \delta_n \frac{e^{lt} - 1}{l} + |(I-P_n) x_m(t)| \leq
     \delta_n \frac{e^{lT} - 1}{l} + |(I-P_n)W|
\end{eqnarray*}
This shows that $\{x_n\}$ is a Cauchy sequence in ${\cal
C}([0,T],H)$, hence it converges uniformly to $x^*:[0,T]\to W$.
From lemma \ref{lem:limit} it follows that  $\frac{d
x^*}{dx}=F(x)$.

{\bf Uniqueness.} Let $x:[0,T] \to W$ be a solution of
(\ref{eq:feqLN}) with the initial condition $x(0)=x_0$. We will
show that $x_n $ converge to $x$. We apply lemma
\ref{lem:estmLogN} to $n$-dimensional projection and the function
$P_n x(t)$. We obtain
\begin{equation}
  |x_n(t) - P_n(x(t))|_n \leq  \delta_n \frac{e^{lt} - 1}{l}.
\end{equation}
Since the tail $(I-P_n)x(t)$ is uniformly bounded we see that $x_n
\to x$ uniformly.

{\bf Lipschitz constant on $W$}. From equation
(\ref{eq:estmdiffg}) applied to $m=n$ for different initial
conditions (we denote the functions by $x_n$ and $y_n$ and the
initial conditions $x_0$ and $y_0$) we obtain
\begin{equation}
  |x_n(t) - y_n(t)| \leq e^{lt}|P_n x_0 - P_n y_0| + \delta_n \frac{e^{lt} - 1}{l}
  \label{eq:diffnlip}
\end{equation}
Let $x_n \to x$ and $y_n \to y$. Then passing to the limit in
(\ref{eq:diffnlip}) gives
\begin{equation}
   |x(t) - y(t)| \leq e^{lt}|x_0 - y_0|.
\end{equation}

Assertion 4 follows easily from  from the previous ones. \qed

%% file: exproof.tex
\section{Existence theorems for Navier-Stokes system in 2D}
\subsection{Some easy lemmas about Fourier series}
\begin{lemma}
\label{lem:Fourier1}
  Let $u \in C^n({\Bbb T}^d,{\Bbb C})$ and let $u_k$ for $k \in {\Bbb Z}^d$ be
  a Fourier coefficient of $u$. Then
  there exists $M$
  \begin{displaymath}
     |u_k| \leq \frac{M}{|k|^n}
  \end{displaymath}
\end{lemma}

\begin{lemma}
\label{lem:Fourier2}
  Assume that $|u_k| \leq \frac{M}{|k|^{\gamma}}$ for $k \in {\Bbb Z}^d$. Let  $n \in {\Bbb N}$
  be  such that $\gamma - n > d$,
  then the function $u(x)=\sum_{k \in {\Bbb Z}^d} u_k e^{ikx}$ belongs to $C^n({\Bbb T}^d,{\Bbb C})$.
  The series
  \begin{displaymath}
    \frac{\partial^s u}{\partial x_{i_1} \dots x_{i_s}}=
    \sum_{k \in {\Bbb Z}^d} u_k  \frac{\partial^s }{\partial x_{i_1} \dots x_{i_s}} e^{ikx}
  \end{displaymath}
  converge uniformly for $0 \leq s \leq n$.
\end{lemma}

\begin{lemma}
\label{lem:Fourier3}
  Assume that for some $\gamma>0$, $a>0$ and $D>0$ we have
  $|u_k| \leq \frac{D e^{-a|k|}}{|k|^{\gamma}}$ for
  $k \in {\Bbb Z}^d \setminus \{0\}$.

  Then the function $u(x)=\sum_{k \in {\Bbb Z}^d} u_k e^{ikx}$ is analytic.
\end{lemma}

Let $H=\left\{ \{ u_k\}\ |  \quad \sum_{k \in {\Bbb Z}^d} |u_k|^2
< \infty \right\}$. Obviously $H$ is a Hilbert space. Let $F$ be
the right side of (\ref{eq:NSgal1})
\begin{equation}
 F(u)_k= -i \sum_{k_1}(u_{k_1}|k)\sqcap_k u_{k-k_1} - \nu k^2u_k + \sqcap_k f_k
\end{equation}
For a general $u \in H$ we cannot claim that $F(u) \in H$. But
when $|u_k|$ decreases fast enough we have the following
\begin{lemma}
\label{lem:compt}
 Let $W(D,\gamma)=\left\{ u \in H \ | \ |u_k| \leq \frac{D}{|k|^\gamma}
\right\}$, then
\begin{description}
  \item[1.] if $\gamma > \frac{d}{2}$, then $W(D,\gamma)$ satisfies condition C2.
  \item[2.] if $\gamma - 2> \frac{d}{2}$ and $\gamma > d$, then the
   function $F:W(D,\gamma)\to H$ is continuous and condition C3 is satisfied on
     $W(D,\gamma)$.
  \item[3.] if $\gamma >d+1$, then  condition D is satisfied on $W(D,\gamma)$.
\end{description}
\end{lemma}
{\bf Proof:} To prove assertion 1 it is enough to show that
$W(d,\gamma)$ is bounded, closed (obvious) and is is componentwise
bounded by some $v=\{v_k\}$, such that $v \in H$. We set
$v_k=\frac{D}{|k|^\gamma}$. Observe that $v \in H$, because
\begin{equation}
  \sum_{k \in {\Bbb Z}^d} |v_k|^2 \leq C D^2 \sum_{n=1}^\infty \frac{n^{d-1}}{n^{2\gamma}}
\end{equation}
and the series converge when $2\gamma - (d-1) > 1$. This concludes
the proof of assertion 1.

To prove assertion 2 we can assume that $f=0$ (it is just a
constant vector in $H$). From lemma \ref{lem:estmQ} if follows
immediately that for $u \in W$ we obtain
\begin{displaymath}
|F(u)_k | \leq \frac{C}{|k|^{\gamma-1}} + \frac{\nu
D}{|k|^{\gamma-2}} \leq \frac{B}{|k|^{\gamma-2}}.
\end{displaymath}
Hence $F(u) \in W(B,\gamma-2) \subset H$, when $\gamma-2 >
\frac{d}{2}$. Hence $F(W(D,\gamma)) \subset W(B,\gamma-2)$. Since
the convergence in $W(B,\gamma-2)$ is equivalent to componentwise
convergence, the same holds for the continuity. It is obvious that
$F(u)_k$ continuous on $W(d,\gamma)$, because the series defining
it is uniformly convergent, hence $F$ is continuous  on
$W(d,\gamma)$.

We prove now assertion 3. Observe that
\begin{equation}
  \frac{\partial N_k}{\partial u_{k_1}}=(\cdot| k)\sqcap_k u_{k - k_1} +
      (u_{k-k_1}|k)\sqcap_k
\end{equation}
We will treat here $u_k$ as one dimensional object, but the
argument is generally correct, i.e. treating $u_k$ as a vector
will introduce only an additional constant and will not affect the
proof. We estimate
\begin{equation}
  \left|\frac{\partial N_k}{\partial u_{k_1}}\right|(W) \leq \frac{2D|k|}{|k-k_1|^\gamma}
\end{equation}

Hence the sum, $S(k)$, appearing in condition D can be estimated
as follows
\begin{eqnarray*}
 S(k)=1/2\sum_{k_1 \in {\Bbb Z}^d \setminus \{0,k\}} {\left|\frac{\partial N_k}{\partial u_{k_1}}\right|}(W) +
  1/2 \sum_{k_1 \in {\Bbb Z}^d \setminus \{0,k\}} {\left|\frac{\partial N_{k_1}}{\partial u_k}\right|}(W) \leq \\
    D|k| \sum_{k_1 \in {\Bbb Z}^d \setminus \{0,k\}} \frac{1}{|k-k_1|^\gamma}  +
    D \sum_{k_1 \in {\Bbb Z}^d \setminus \{0,k\}} \frac{|k_1|}{|k-k_1|^\gamma}
\end{eqnarray*}
Now observe that
\begin{equation}
\sum_{k_1 \in {\Bbb Z}^d \setminus \{0,k\}}
\frac{1}{|k-k_1|^\gamma} < \sum_{k_1 \in {\Bbb Z}^d, k_1 \neq 0}
\frac{1}{|k|^\gamma}=C(d,\gamma) < \infty, \quad \mbox{ for
$\gamma > d$}
\end{equation}

To estimate the sum $\sum_{k_1 \in {\Bbb Z}^d \setminus \{0,k\}}
\frac{|k_1|}{|k-k_1|^\gamma}$ we show that there exists a constant
$A$, such that
\begin{equation}
\frac{|k_1|}{|k-k_1|} < A|k|, \quad \mbox{for $k,k_1 \in {\Bbb
Z}^d \setminus \{0\} $, $k \neq k_1$}.
\end{equation}
Observe that for $|k_1| \leq 2|k|$, $k_1 \neq 0$, $k_1 \neq k$ we
can estimate the denominator by $1$  hence we have
\begin{equation}
\frac{|k_1|}{|k-k_1|} \leq 2|k|.
\end{equation}
For $|k_1| > 2 |k|$ we have
\begin{equation}
  \frac{|k_1|}{|k-k_1|}=\frac{1}{\left|\frac{k_1}{|k_1|} - \frac{k}{|k_1|} \right|} \leq
     \frac{1}{1 -  \frac{|k|}{|k_1|} } \leq 2.
\end{equation}
So we can take $A=2$.

Now we  make the following estimate
\begin{equation}
 \sum_{k_1 \in {\Bbb Z}^d \setminus \{0,k\}} \frac{|k_1|}{|k-k_1|^\gamma} \leq
   A|k| \sum_{k_1 \in {\Bbb Z}^d \setminus \{0,k\}} \frac{1}{|k-k_1|^{\gamma-1}}
   < A C(d,\gamma-1) |k|,
\end{equation}
provided $\gamma-1 > d$.

So we have $S(k) < \left(D C(d,\gamma) + A D
C(d,\gamma-1)\right)|k|$ and since $\lambda_k=-\nu |k|^2$, we see
that $l$ satisfying condition {\em D} exists. \qed

\subsection{Existence theorems}
We set the dimension $d=2$. We again assume that the force $f$ is
such that $f_k=0$ for $|k| > K$ (in \cite{MS}) more general force
is treated.

Observe that from lemma \ref{lem:compt} it follows that to have
conditions C1, C2, C3, D on the trapping regions constructed in
section \ref{sec:trap} we need $\gamma > 3$.

\begin{theorem}
  \label{thm:exist1}
  If for some $D$ and $\gamma >3$
  \begin{equation}
    |u_k(0)| \leq \frac{D}{|k|^\gamma}
  \end{equation}
  then the solution of (\ref{eq:NSgal1}) is defined for all $t>0$ and
  there exists a constant $D'$, such that
  \begin{equation}
    |u_k(t)| \leq \frac{D'}{|k|^\gamma}, \quad t>0.
  \end{equation}
\end{theorem}

The following theorem tells that if we start with analytic initial
conditions that the solution will remain analytic (in space
variables).
\begin{theorem}
\label{thm:exist2}
  If for some $D$, $\gamma >3$ and $a >0$
  \begin{equation}
    |u_k(0)| \leq \frac{D}{|k|^\gamma}e^{-a|k|}
  \end{equation}
  then the solution of (\ref{eq:NSgal1}) is defined for all $t>0$ and
  there exist  constants $D'$ and $a'>0$  such that
  \begin{equation}
    |u_k(t)| \leq \frac{D'}{|k|^\gamma}e^{-a'|k|}, \quad t>0
  \end{equation}
\end{theorem}

Next theorem states that the solution starting from regular
initial conditions becomes analytic immediately.
\begin{theorem}
\label{thm:exist3}
  Assume that for some $D$, $\gamma >3$ and $a >0$ the initial conditions satisfy
  \begin{equation}
    |u_k(0)| \leq \frac{D}{|k|^\gamma}
  \end{equation}
  then the solution of (\ref{eq:NSgal1}) is defined for all $t>0$ and
  and for any  $t_0 >0$ one can find  constants $D'$ and $a'>0$  such that
  \begin{equation}
    |u_k(t_0)| \leq \frac{D'}{|k|^\gamma}e^{-a'|k|}
  \end{equation}
\end{theorem}

\noindent {\bf Proof of theorem \ref{thm:exist1}:} Observe first
that the enstrophy of $\{u_k(0)\}$ is finite. Let take $V_0
> \max(V(\{u_k\}), V^*)$. From theorem \ref{thm:trap1} it follows
that there exists $K$ and $D'$, such that $\{u_k(0)\}$ belongs to
the trapping set $N=N(V_0,K,\gamma,D')$. Observe that $N \subset
W(D',\gamma)$, hence we can pass to the limit with solutions
obtained from Galerkin projections (see theorem
\ref{thm:limitLN}). \qed

\noindent {\bf Proof of theorem \ref{thm:exist2}:} The proof is
essentially the same as for theorem \ref{thm:exist1}, the only
difference is: we use theorem \ref{thm:trap2} instead of theorem
\ref{thm:trap1}. \qed

\noindent {\bf Proof of theorem \ref{thm:exist3}:} The global
existence was proved in theorem \ref{thm:exist1}. To prove the
estimate for $|u_k(t_0)|$ we use theorem \ref{thm:trap3} to obtain
\begin{equation}
  |u_k(t_0)| \leq \frac{D'}{|k|^\gamma}e^{-a|k|t_0},
\end{equation}
which finishes the proof. \qed

\begin{theorem}
$d=2$.  If $u_0 \in C^5$ then the classical solution of NS
equations such that $u(0,x)=u_0(x)$ exists for all $t>0$ and it is
analytic in space variables for $t>0$.
\end{theorem}
{\bf Proof:} From lemma \ref{lem:Fourier1} it follows that the
Fourier coefficients of $u_0$, $\{u_{0,k}\}$, satisfy assumptions
of the theorem \ref{thm:exist1} with $\gamma=5$. Hence there
exists a solution, $\{u_k(t)\}$, of (\ref{eq:NSgal1}) in $H$, such
that $u_k(0)=u_{0,k}$.

Let us set $u(t,x)=\sum_{k \in {\Bbb Z}^2 \setminus \{0\}} u_k(t)
e^{ikx}$. It is easy to see that $u(t,x)$ is a classical solution
of the Navier-Stokes system, because the Fourier series for all
terms in the NS equations converge fast enough (compare proof of
lemma \ref{lem:ladder}).

From the theorem \ref{thm:exist3} and lemma \ref{lem:Fourier3} it
follows that the function $u(t_0,\cdot)$ is analytic in space
variables for any $t_0>0$. \qed

The following theorem is an easy consequence of theorem
\ref{thm:limitLN}.
\begin{theorem}
Assume $d=2$ and $\gamma>3$. Let $W$ be any of the trapping
regions defined in theorems \ref{thm:trap1} and \ref{thm:trap2},
then the Navier-Stokes system induces a semidynamical system on
$W$.
\end{theorem}

%% file: sankal.tex
\section{Trapping regions  in 3D}
\label{sec:SanKal}

The goal of this section is to present method by Sannikov and
Kaloshin \cite{S} for constructing a trapping region  for small
initial data.

Let us state a result, which is not contained in \cite{S}, but can
be easily obtained using the technique presented there.

We set the dimension $d=3$. We  assume the force $f$ is zero.
\begin{theorem}
  \label{thm:exist3d}
  For any $\gamma > 3.5$, there exists $D_0=D_0(\gamma,\nu)$
  such that for all $D < D_0$, if
  \begin{equation}
    |u_k(0)| \leq \frac{D}{|k|^\gamma}
  \end{equation}
  then the solution of (\ref{eq:NSgal1}) is defined for all $t>0$ and
  \begin{equation}
    |u_k(t)| \leq \frac{D}{|k|^\gamma}, \quad t>0
  \end{equation}
\end{theorem}
{\bf Proof:} Let
\begin{equation}
W=\left\{ \{u_k\} \: | \:  |u_k| \leq \frac{D}{|k|^\gamma}
\right\}.
\end{equation}

From lemma \ref{lem:estmQ} it follows that for $\{u_k\} \in W$ we
have
\begin{equation}
   \frac{d |u_k|}{dt} \leq \left|\sum (u_{k_1}|k) \sqcap_k u_{k-k_1}\right| - \nu |k|^2 |u_k|
    \leq \\
   \frac{D^2 C_Q(3,\gamma)}{|k|^{\gamma -1}} - \nu |k|^2 |u_k|.
\end{equation}
Hence $W$ is a trapping region if for every $k$ we have
\begin{equation}
   \frac{D^2 C_Q(3,\gamma)}{|k|^{\gamma -1}} -  \frac{\nu D}{|k|^{\gamma-2}} < 0.
\end{equation}
We obtain
\begin{equation}
  \frac{D C_Q(3,\gamma)}{\nu} < |k|, \quad k \in {\Bbb Z}^3 \setminus \{0\}.
\end{equation}
Hence if
\begin{equation}
  D < D_0=\frac{\nu}{C_Q(3,\gamma)},
\end{equation}
then $W$ is a trapping region for all projections of the Navier
Stokes equations. From lemma \ref{lem:compt} it follows that the
conditions C1,C2,C3 are satisfied (it is easy to see that
condition D holds if $\gamma>4$.) Hence we can pass to the limit
with the dimension of Galerkin projection to obtain a desired
solution. \qed

One can easily state similar theorem for analytic initial
condition.

Let us comment on the Sannikov and Kaloshin result \cite{S}. They
constructed the trapping region of the form $|u_k| \leq
\frac{D}|k|^2 e^{-v|k|t}$, $t\geq 0$, where $v>0$. The methods
developed in this paper require more compactness  at $t=0$ to be
directly applicable to this trapping region.